# Are there parts of our arithmetical competence that no sound formal system can duplicate?

*Bhupinder Singh Anand*

*In a review of Roger Penrose's "Shadows of the Mind", David Chalmers opined as implausible that there could be parts of our arithmetical competence that no sound formal system could ever duplicate. We prove, however, that the recursive number-theoretic relation $x=Sb(y\ 19|Z(y))$ - which is algorithmically verifiable since the recursive function $Sb(y\ 19|Z(y))$ is Turing-computable - cannot be introduced through definition in any consistent formal system of Arithmetic.*

**Introduction**

In a review of Roger Penrose's "Shadows of the Mind" [Pe94], David Chalmers remarked ([Ch95], §2.5):

> "... I have some sympathy with Penrose's idea that we have an underlying sound competence, even if our performance sometimes goes astray. But further, it seems to me that to hold that this is the only problem in Penrose's argument would be to concede too much power to the argument. It would follow, for example, that there are parts of our arithmetical competence that no sound formal system could ever duplicate; it would seem that our unsoundness would be essential to our capacity to see the truth of Gödel sentences, for example. This would be a remarkably strong conclusion, and does not seem at all plausible to me".

The question arises: Is there a number-theoretic relation that is algorithmically verifiable - and can thus be viewed as "part of our arithmetical competence" - but which cannot be consistently "duplicated" within any "sound formal system" of Arithmetic?

We address this question, somewhat obliquely, in this paper, and briefly note some immediate consequences. We essentially prove that the recursive[1] number-theoretic relation $x=Sb(y\ 19|Z(y))$[2] - which is algorithmically verifiable since the recursive function $Sb(y\ 19|Z(y))$ is Turing-computable[3] - cannot be introduced through definition in any consistent[4] formal system of Arithmetic such as Gödel's formal system P [Go31a], or standard PA[5].

We note that, in his seminal 1931 paper [Go31a], Gödel defines the recursive number-theoretic function $Sb(x\ 19|Z(y))$ as the Gödel-number[6] of the formula[7], of his formal system P, that is obtained from the P-formula whose Gödel-number is $x$ by substituting the numeral[8] $[y]$[9], whose Gödel-number is $Z(y)$, for the variable whose Gödel-number is 19 wherever the latter occurs free in the P-formula whose Gödel-number is $x$.

---

[1] We follow Gödel's definition of recursive number-theoretic functions and relations ([Go31a], p14-17).

[2] This relation occurs on 2nd August 2002 in the correspondence titled "The Godel's Loop" between Antonio Espejo <CASAFARFARA@terra.es> and "Rupert" <rupertmccallum@yahoo.com> in the Google newsgroup sci.logic.

[3] We call a function $F$ "Turing-computable" if, and only if, there is an algorithm $U$ that computes $F$ as defined by Mendelson ([Me64], p231).

[4] We take the classical definition of "consistency" as defined by Mendelson ([Me64], p37).

[5] We take standard PA to be the first order theory S defined by Mendelson ([Me64], p102), in which addition and multiplication are the standard interpretations of the primitive symbols "+" and "*" respectively.

[6] By the "Gödel-number" of a formula of P, we mean the natural number corresponding to the formula in the 1-1 correspondence defined by Gödel ([Go31a], p13).

[7] By "formula", we shall henceforth mean a "well-formed formula" as defined by Gödel ([Go31a], p11).

[8] We follow Gödel's definition of a numeral ([Go31a], p10).

## 1 A Meta-theorem of recursive asymmetry

*Definition*: A recursive number-theoretic function or relation is asymmetrical in P if it is not the standard interpretation[10] of any of its formal representations[11] in P.

*Meta-theorem 1*: There is a recursive number-theoretic relation that is asymmetrical in P.

*Proof*: We consider the number-theoretic relation $x=Sb(y\ 19|Z(y))$.

(*a*) We assume that no recursive number-theoretic function or relation is asymmetrical in P. In other words, we assume that every recursive number-theoretic function or relation is the standard interpretation of at least one of its formal representations in P.

(*b*) Let the P-formula $[F(x, y)]$ denote a formal representation of the recursive number-theoretic relation $x=Sb(y\ 19|Z(y))$.

(*c*) Let $F(x, y)$ denote the standard interpretation[12] of $[F(x, y)]$.

(*d*) We consider the case where $x=Sb(y\ 19|Z(y))$ is an abbreviation[13] for $F(x, y)$. In other words, we assume that, if we use Gödel's recursive definitions ([Go31a],

---

[9] Henceforth, we use square brackets to differentiate between a formula $[R(x)]$ and its interpretation $R(x)$ under some given interpretation, say M, where we follow Mendelson's definition of "interpretation" ([Me64], p48, §2). We note that the interpreted relation $R(x)$ is obtained from the formula $[R(x)]$ by replacing every primitive, undefined symbol of P in the formula $[R(x)]$ by an interpreted mathematical symbol under the interpretation M (i.e. a symbol that is a shorthand notation for some, semantically well-defined, concept of classical mathematics under the interpretation M). So, assuming the given interpretation M follows standard logical terminology, the P-formula $[(Ax)R(x)]$ would interpret as the M-proposition $(Ax)R(x)$, and the P-formula $[\sim(Ax)R(x)]$ as the M-proposition $\sim(Ax)R(x)$.

[10] We follow Mendelson's definition of "standard interpretation" ([Me64], p107).

[11] We follow Mendelson's definitions of "expressibility" ([Me64], p117), and "representability" ([Me64], p118), of number-theoretic functions and relations in a formal system such as P.

[12] By definition, $F(x, y)$ is some interpretation of the P-formula $[F(x, y)]$.

[13] Cf. ([Go31a], p11, footnote 22); ([Me64], p31, footnote 1); ([Me64], p82).

p17-20), and follow the reasoning he outlines in Theorem V ([Go31a], p23), we can transform the relation $x=Sb(y\ 19|Z(y))$ into a relation $F(x, y)$, such that all the symbols that occur in $F(x, y)$ are standard interpretations of primitive symbols of P.

Now:

(*i*)  Let $k$ be the Gödel-number of $[F(x, y)]$.

(*ii*)  Then $Sb(k\ 19|Z(k))$ is the Gödel-number of the P-formula $[F(x, k)]$ that we get when, in the P-formula $[F(x, y)]$, we replace the variable $[y]$, wherever it occurs, by the numeral $[k]$.

(We assume that $[y]$ is Gödel-numbered as $19^{14}$ in the Gödel-numbering that yields $k$ as the Gödel-number of $[F(x, y)]$)

(*iii*)  Let $l= Sb(k\ 19|Z(k))$.

(*iv*)  By definition, the instantiation of the number-theoretic function $Sb(y\ 19|Z(y))$, whose abbreviation is $Sb(k\ 19|Z(k))$, must contain an explicit bound $k'$ that is equal to, or larger than, $l$.

This follows from the constructive definition of Gödel's recursive functions ([Go31a], p17, footnote 34). Thus $Sb(y\ 19|Z(y))$ is of the form "$(ix)((x =< f_1(y))$ & $g_1(x, y))$", where "$f_1(y)$" and "$g_1(x, y)$" are abbreviations, for a recursive function and a recursive relation, respectively, of lower rank ([Go31a], p16, Theorem IV). Similarly, $f_1(y)$ is of the form "$(ix)((x =< f_2(y))$ & $g_2(x, y))$", etc.

---

[14] Following the convention set by Gödel in Theorems V and VI (cf. [Go31a], p22, footnote 38), we assign 17 as the Gödel-number of "$x$", and 19 as the Gödel-number of "$y$".

5(Here, "$(ix)R(x)$" denotes the smallest natural number for which the relation $R(x)$ holds, and "$=<$" denotes the relation of "equal to or less than".)

We thus have a finite sequence $f_1(y), f_2(y), ..., f_n(y)$ of recursive functions, of decreasing rank, such that:

$$Sb(y\ 19|Z(y)) =< f_1(y) =< f_2(y) =< ... =< f_n(y),$$

where $n$ is less than or equal to the rank $r$ of $Sb(y\ 19|Z(y))$, and where $f_n(y)$ does not contain any abbreviations.

It follows that $f_n(y)$ occurs, as a bound, in the unabbreviated number-theoretic function whose abbreviation is $Sb(y\ 19|Z(y))$. By Gödel's reasoning, as outlined in Theorem V, $f_n(y)$ is, therefore, either a constant, or of the form "$y+q$" ([Go31a], p23), where $q$ is a natural number that depends on the rank $r$ of $Sb(y\ 19|Z(y))$.

(We note that the representation of $f_n(y)$ in P, as envisaged in Theorem V, would be Gödel's term of the first type ([Go31a], p10), denoted by $[S^q y]$, where '$S^q$' denotes the pre-fixing of the primitive (*successor*) symbol $[S]$[15] of P, to $[y]$, $q$ times.)

We thus have that:

$$l = Sb(k\ 19|Z(k)) =< f_n(k) = k'.$$

(*v*)  Since the standard interpretation of $[F(x, k)]$ is $F(x, k)$, it follows, from (*d*), that $F(x, k)$ is the number-theoretic relation whose abbreviation is $x = Sb(k\ 19|Z(k))$.

---

[15] We note that Gödel uses the symbol "*f*" instead of "*S*" in his paper [Go31a].

(*vi*) By (*iv*), [$F(x, k)$] must, therefore, contain a bound [$k'$], which interprets as a natural number $k'$ that is larger than $l$.

This is impossible, since a formula cannot contain a numeral that, under interpretation, yields a natural number that is equal to, or larger than, the Gödel-number of the formula[16]. It follows that assumption (*a*) does not hold; this proves the meta-theorem.

## 2 Some consequences

(*a*) It now follows that:

> *Meta-lemma 1*: We cannot introduce a finite number of arbitrary recursive number-theoretic functions and relations as primitive symbols into P without risking inconsistency.
>
> *Proof*: Adding "Sb", and a finite number of other functions and relations in terms of which it is defined ([Go31a], p17-19, Defs. 1-31), as new primitive symbols of P, along with associated defining axioms (cf. [Me64], p82 , §9), would yield a formal system P' in which [$x=Sb(y\ 19|Z(y))$] is a P'-formula. By *Meta-theorem 1*, P' would be inconsistent.

(*b*) It further follows that:

---

[16] This is easily proved since, following, for instance, Gödel's assignment of the natural numbers 1 and 3 to the primitive P-symbols [0] and [$S$] respectively ([Go31a], p13), we have that:

(*i*) the Gödel number $p_1^3\ p_2^3\ ...\ p_q^3\ p_{q+1}$ of the numeral [$S^q 0$], which represents the natural number $q$ in P, is greater than $q$ for all $q >= 0$, where $p_i$ is the $i$'th prime;

(*ii*) for any P-formulas [$F$] and [$G$], the Gödel number of the concatenated P-formula [$FG$] is always greater than, or equal to, the individual Gödel numbers of the P-formulas [$F$] and [$G$].

*Meta-lemma 2*: Although a recursive number-theoretic relation, and the standard interpretation of its formal representation, are always equivalent in their instantiations, they are not always formally equivalent.

*Proof*: Let $F(x_1, ..., x_n)$ be any recursive number-theoretic relation, and $G(x_1, ..., x_n)$ be the standard interpretation of one of its formal representations in P.

We thus have that the two relations are equivalent in their instantiations since, by definition, for any given sequence $<a_1, ..., a_n>$ of natural numbers:

$F(a_1, ..., a_n)$ holds if and only if $G(a_1, ..., a_n)$ holds.

However, by *Meta-lemma 1*, it follows that $[F(x_1, ..., x_n)]$ is not necessarily a P-formula. Hence, we cannot conclude that:

$[F(x_1, ..., x_n) <=> G(x_1, ..., x_n)]$ is a P-formula.

This proves the lemma[17].

(*c*) Can every recursive number-theoretic function or relation be defined constructively within P? The classical argument is expressed by Mendelson ([Me64], p82, §9):

"In mathematics, once we have proved, for any $y_1, ..., y_n$, the existence of a unique object $u$ having the property $A(u, y_1, ..., y_n)$, we often introduce a new function $f(y_1, ..., y_n)$ such that $A(f(y_1, ..., y_n), y_1, ..., y_n)$ holds for all $y_1, ..., y_n$. ... It is generally acknowledged that such definitions, though convenient, add nothing really new to the theory."

---

[17] We note that Gödel asserts at the end of his Theorem VII that the equivalence in *Meta-lemma 2* can always be formulated within the formal system P ([Go31a], p31).

More precisely, he argues, in his Proposition 2.29 ([Me64], p82), that, classically, abbreviations for *strongly* representable ([Me64], p118) number-theoretic functions may be introduced as primitive symbols into P, since they can always be eliminated. It follows that:

> *Meta-lemma 3*: If every strongly representable number-theoretic function can be introduced as a new function letter into P, without affecting the consistency of P, then every recursive function is not *strongly* representable[18] in P.

(*d*) We also note that every arithmetical relation is the standard interpretation of one of its representations in a formal system of Arithmetic such as standard PA. Since *Meta-Lemma 1* can be seen to hold in any such system, it follows from *Meta-lemma 2* that:

> *Meta-lemma 4*: Although every primitive recursive relation is equivalent to some arithmetical relation, in the sense that they are always equivalent in their instantiations, such equivalence cannot be formulated within a formal system of Arithmetic such as standard PA.

---

[18] By Mendelson's reasoning in his Proposition 3.23 ([Me64], p131), every recursive function is representable in P.

(*Acknowledgements: My thanks to correspondents in various newsgroups who offered comments and counter-arguments that eventually led to the arguments of Meta-theorem . My grateful thanks to Dr. Damjan Bojadziev for comments on the formatting and presentation. Author's e-mail: anandb@vsnl.com.*)

(*Updated: Saturday 10$^{th}$ May 2003 1:23:28 AM IST by re@alixcomsi.com*)